\documentclass{amsart}
\usepackage{amsmath, amssymb, latexsym, graphicx}
\title{Approximations to the Kruskal-Katona Theorem}
\author{Andrew Frohmader}
\address{Department of Mathematics, 581 Malott Hall, Cornell University, Ithaca, NY 14853-4201}
\email{froh@math.cornell.edu}

\newtheorem{theorem}{Theorem}[section]

\newtheorem{corollary}[theorem]{Corollary}
\newtheorem{lemma}[theorem]{Lemma}
\newtheorem{definition}[theorem]{Definition}

\def\proof{\smallskip\noindent {\it Proof: \ }}
\def\endproof{\hfill\ensuremath{\square}\medskip}

\begin{document}

\maketitle

\begin{abstract}
Approximations to the Kruskal-Katona theorem are stated and proven.  These approximations are weaker than the theorem, but much easier to work with numerically.
\end{abstract}

\section{Introduction}

The Kruskal-Katona theorem \cite{kruskal,katona} characterizes the face vectors of simplicial complexes.  It gives sharp bounds, and given a proposed face vector, can either produce a complex with that exact face vector or else say that no such complex exists.

Unfortunately, the bounds are somewhat awkward to work with.  While a computer can readily be programmed to give the precise bounds, they aren't that numerically intuitive.  In this paper, we give somewhat weaker bounds that are much easier to work with.  We state the bounds now, and return to prove them later.

\begin{theorem}  \label{withoutr}
Let $\Delta$ be a simplicial complex and let $k > p > 0$ be integers.  Suppose further that $f_{k-1}(\Delta) > 0$.  Then
$$f_{p-1}(\Delta) > {(k!)^{{p \over k}} \over p!} \bigg(1 + {k-p \over 2\sqrt[k]{k! f_{k-1}(\Delta)}} \bigg)^p (f_{k-1}(\Delta))^{{p \over k}}.$$
\end{theorem}

\begin{corollary}  \label{noreasy}
Let $\Delta$ be a simplicial complex and let $k > p > 0$ be integers.  Suppose further that $f_{p-1}(\Delta) > 0$.  Then
$$f_{p-1}(\Delta) > {(k!)^{{p \over k}} \over p!} (f_{k-1}(\Delta))^{{p \over k}}.$$
\end{corollary}

\begin{corollary}  \label{symmetric}
Let $\Delta$ be a simplicial complex with $k > 0$ and $f_{k-1}(\Delta) > 0$.  Then
$$f_0(\Delta) > \sqrt{2! f_1(\Delta)} > \sqrt[3]{3! f_2(\Delta)} > \sqrt[4]{4! f_3(\Delta)} > \dots > \sqrt[k+1]{(k+1)! f_k(\Delta)}.$$
\end{corollary}

\begin{theorem}  \label{colorapprox}
Let $r \geq k > p > 0$ be integers and let $\Delta$ be an $r$-colorable simplicial complex.  Then
$$f_{p-1}(\Delta) \geq {r \choose p}{r \choose k}^{-{p \over k}} (f_{k-1}(\Delta))^{{p \over k}}.$$
\end{theorem}

Theorem~\ref{colorapprox} is equivalent to \cite[Theorem 5.1]{balanced}.  We give a proof that is independent of theirs and simpler.

\begin{theorem}  \label{withr}
Let $\Delta$ be a simplicial complex and let $r \geq k > p > 0$ be integers.  If $f_{k-1}(\Delta) \leq {r \choose k} + {r-1 \choose k-1}$, then
$$f_{p-1}(\Delta) \geq {r \choose p}{r \choose k}^{-{p \over k}} (f_{k-1}(\Delta))^{{p \over k}}.$$
\end{theorem}

\begin{corollary}  \label{flagdim}
Let $\Delta$ be a flag complex of dimension $r-1$.  Then
$$f_{p-1}(\Delta) \geq {r \choose p}{r \choose k}^{-{p \over k}} (f_{k-1}(\Delta))^{{p \over k}}.$$
\end{corollary}

\begin{corollary}  \label{flagnodim}
Let $\Delta$ be a flag complex and let $r \geq k > p > 0$ be integers.  If $f_{k-1}(\Delta) < {r+1 \choose k}$, then
$$f_{p-1}(\Delta) \geq {r \choose p}{r \choose k}^{-{p \over k}} (f_{k-1}(\Delta))^{{p \over k}}.$$
\end{corollary}

This paper is structured as follows.  In Section~2, we introduce the Kruskal-Katona theorem and an approximation to it due to Lov\'{a}sz.  We explain why these are numerically awkward to work with, why Theorem~\ref{withoutr} and Corollaries~\ref{noreasy} and~\ref{symmetric} are much easier, and prove these bounds.  In Section~3, we introduce the notion of colored complexes and the Frankl-F\"{u}redi-Kalai theorem, and then use it to prove Theorems~\ref{colorapprox} and~\ref{withr}.  In Section~4, we introduce the notion of a flag complex and prove Corollaries~\ref{flagdim} and~\ref{flagnodim}.  In Section~5, we give some graphs that numerically compare the various bounds discussed throughout the paper.

\section{Avoiding binomial coefficients}

In this section, we give some background material on the problem of characterizing the face vectors of simplicial complexes.  The Kruskal-Katona theorem, Theorem~\ref{kktheorem}, solves this problem.  Theorem~\ref{continuouskk} gives a numerical approximation to the Kruskal-Katona theorem.  We explain why these can be awkward and numerically unintuitive, while Theorem~\ref{withoutr} and Corollaries~\ref{noreasy} and~\ref{symmetric} are much easier to work with.  We also prove Theorem~\ref{withoutr} and Corollaries~\ref{noreasy} and~\ref{symmetric}.

Recall that a \textit{simplicial complex} $\Delta$ on a vertex set $W$ is a collection of subsets of $W$ such that (i) for every $v \in W$, $\{v\} \in \Delta$ and (ii) for every $B \in \Delta$, if $A \subset B$, then $A \in \Delta$.  The elements of $\Delta$ are called \textit{faces}.  A face on $i$ vertices is said to have \textit{dimension} $i-1$, while the dimension of a complex is the maximum dimension of a face of the complex. A face of the same dimension as the complex is called a \textit{facet}.  If all faces that are maximal with respect to inclusion are facets, we say that the complex is \textit{pure}.

The \textit{$i$-th face number} of a simplicial complex $\Delta$, $f_{i-1}(\Delta)$, is the number of faces of $\Delta$ on $i$ vertices.  The \textit{face vector} $f(\Delta)$ of $\Delta$ lists the face numbers of $\Delta$.

One can ask which integer vectors can arise as face vectors of simplicial complexes.  The Kruskal-Katona theorem \cite{kruskal,katona} answers this question, and gives a complete characterization of face vectors of simplicial complexes.  The statement of the theorem requires a lemma.

\begin{lemma} \label{kklemma}
Given positive integers $m$ and $k$, there is a unique way to pick integers $s \geq 0$ and $n_k, n_{k-1}, \dots, n_{k-s}$ such that
$$m = {n_k \choose k} + {n_{k-1} \choose k-1} + \dots + {n_{k-s} \choose k-s}$$
and $n_k > n_{k-1} > \dots > n_{k-s} \geq k-s > 0$.
\end{lemma}

\begin{theorem}[Kruskal-Katona] \label{kktheorem}
Let $\Delta$ be a simplicial complex and let $m = f_{k-1}(\Delta)$.  Let
$$m = {n_k \choose k} + {n_{k-1} \choose k-1} + \dots + {n_{k-s} \choose k-s}$$
as in Lemma~\ref{kklemma}.  Then
$$f_{k-2}(\Delta) \geq {n_k \choose k-1} + {n_{k-1} \choose k-2} + \dots + {n_{k-s} \choose k-s-1}.$$
Furthermore, if a positive integer vector with 1 as its first entry satisfies these inequalities for all $k \geq 1$, then it is the face vector of some simplicial complex.
\end{theorem}

The constants of Lemma~\ref{kklemma} are practical to compute.  We pick $n_k$ such that ${n_k \choose k} \leq m < {n_k + 1 \choose k}$, and then repeat with $k-1$ instead of $k$ and $m - {n_k \choose k}$ instead of $m$.  We stop when we would end up with 0 as the new value for $m$, and this is guaranteed to happen no later than when we use 1 for $k$.

It is not difficult to show that ${(n - k + 1)^k \over k!} \leq {n \choose k} \leq {(n - .5k + .5)^k \over k!}$ and that ${n \choose k}$ is usually much closer to the upper bound than the lower bound.  We can usually set ${(n - .5k + .5)^k \over k!} = m$ and solve for $n$ to get a pretty good approximation to $n_k$, check a few nearby integer values to get the exact value of $n_k$, then repeat with $n_{k-1}$, and so on.

The bounds of Theorem~\ref{kktheorem} can be chained so that if $i \geq 2$, we get
$$f_{k-i}(\Delta) \geq {n_k \choose k-i+1} + {n_{k-1} \choose k-i} + \dots + {n_{k-s} \choose k-i-s+1}.$$
While not technically the formula we would get by chaining the bounds of the Kruskal-Katona theorem if $k-i-s+1 < 0$, this is numerically equivalent if we follow the convention that ${n \choose k} = 0$ when $k < 0$.

Unfortunately, computing the precise numerical bounds is somewhat involved.  We may have to compute as many as $k$ values of $n_i$, with each one after the first depending rather chaotically on previous computations.  This makes the numerical values rather unintuitive.

There is one effort due to Lov\'{a}sz \cite{lovasz} at simplifying this considerably by using only one binomial coefficient rather than up to $k$ of them.  This approximation doesn't give sharp bounds, but is reasonably close.  We need a definition in order to state the approximation.

\begin{definition}
\textup{Let $x$ be a real number and let $k$ be a positive integer.  Define}
$${x \choose k} = {x(x-1) \dots (x-k+1) \over k!}.$$
\end{definition}

If $x$ is a positive integer, this corresponds to the usual binomial coefficients.  The point of the definition is to interpolate between integers, so that $x$ can be any other real number.  For example, ${3.5 \choose 2} = 4.375$.

It is clear from the definition that ${x \choose k}$ is a polynomial in $x$ of degree $k$. Its zeroes are $x = 0, 1, \dots, k-1$, and $\lim_{x \to \infty} {x \choose k} = \infty$.  In addition, ${x \choose k}$ is strictly increasing on $x > k-1$, as all $k$ of the factors of its numerator are positive and strictly increasing on this domain.  As such, given any real $c > 0$, there is a unique $x > k-1$ such that ${x \choose k} = c$.  This allows us to state the following approximation to the Kruskal-Katona theorem.

\begin{theorem}[Lov\'{a}sz] \label{continuouskk}
Let $\Delta$ be a simplicial complex and let $k > p > 0$ be integers.  If $f_{k-1}(\Delta) = {x \choose k}$ for $x \geq k-1$, then $f_{p-1}(\Delta) \geq {x \choose p}$.
\end{theorem}

This approximation is much easier to state than the full Kruskal-Katona theorem.  It is not sharp in general, and only coincides with the bounds of the Kruskal-Katona theorem when $x$ is an integer.

Unfortunately, it still isn't that easy to compute.  Finding $x$ involves finding the largest real root of a polynomial of degree $k$.  It is not difficult to find a numerical approximation by computer, but it still doesn't give that great of intuition on how large the bounds are otherwise.

We use Theorem~\ref{continuouskk} to prove Theorem~\ref{withoutr}.  First, we need a couple of lemmas.

\begin{lemma}  \label{narrower}
Let $k > p > 0$ be integers and let $x > k-1$ be a real number.  Let $c = {k - p \over 2}$.  Then
$$\sqrt[k]{(x)(x-1) \dots (x-k+1)} < \sqrt[p]{(x-c)(x-c-1) \dots (x-c-p+1)}.$$
\end{lemma}

\proof  Let $d = {k-1 \over 2}$.  First note that if $a > b \geq 0$, then
$$(x-d+a)(x-d-a) = (x-d)^2 - a^2 < (x-d)^2 - b^2 = (x-d+b)(x-d-b).$$
We can take both sides of the inequality in the lemma to the $kp$ power and get that the lemma is equivalent to
$$x^p(x-1)^p \dots (x-k+1)^p < (x-c)^k(x-c-1)^k \dots (x-c-p+1)^k.$$
If we start with the left side, we can attain the right side by a sequence of substitutions that consist of replacing $(x-d+a)(x-d-a)$ by $(x-d+b)(x-d-b)$ with $a > b \geq 0$.  We can do this by at each step choosing $a$ to be the largest value such that we have too many factors of $(x-d+a)(x-d-a)$ and $b$ to be the largest value such that we need more factors of $(x-d+b)(x-d-b)$.  \endproof

\begin{lemma}  \label{subtractc}
Let $p$ be a positive integer and let $x > p-1$ and $c > 0$ be real numbers.  Then
$$\sqrt[p]{(x+c)(x+c-1) \dots (x+c-p+1)} > \sqrt[p]{(x)(x-1) \dots (x-p+1)} + c.$$
\end{lemma}

\proof  We can readily compute
$${d \over dx} \sqrt[p]{(x) \dots (x-p+1)} = {1 \over p} \sqrt[p]{(x) \dots (x-p+1)} \bigg({1 \over x} + \dots + {1 \over x-p+1}\bigg).$$
We apply the geometric mean-harmonic mean inequality to the set $\{x, x-1, \dots, x-p+1\}$ to get
$$\sqrt[p]{(x)(x-1)(x-2) \dots (x-p+1)} > {p \over {1 \over x} + {1 \over x-1} + \dots + {1 \over x-p+1}}.$$
From this, it follows that
$${d \over dx} \sqrt[p]{(x)(x-1)(x-2) \dots (x-p+1)} \geq 1.$$
We use the Fundamental Theorem of Calculus to compute
\begin{eqnarray*}
c & = & \int_x^{x+c} 1\ dt \\ & < & \int_x^{x+c} {d \over dt} \sqrt[p]{(t)(t-1)(t-2) \dots (t-p+1)}\ dt \\ & = & \sqrt[p]{(x+c)(x+c-1) \dots (x+c-p+1)} - \sqrt[p]{(x)(x-1) \dots (x-p+1)},
\end{eqnarray*}
from which the statement of the lemma follows.  \endproof

\smallskip\noindent \textit{Proof of Theorem~\ref{withoutr}: \ }  Let $c = {k - p \over 2}$.  Let $f_{k-1}(\Delta) = {x \choose k}$ as in the statement of Theorem~\ref{continuouskk}.  By the theorem, $f_{p-1}(\Delta) \geq {x \choose p}$.  We can use this to compute
\begin{eqnarray*}
& & f_{p-1}(\Delta) \\ & \geq & {x \choose p} \\ & = & {x \choose p} {x \choose k}^{-{p \over k}} f_{k-1}(\Delta)^{{p \over k}} \\ & = & {(k!)^{{p \over k}} \over p!} {x(x-1) \dots (x-p+1) \over \big( x(x-1) \dots (x-k+1) \big)^{{p \over k}}} f_{k-1}(\Delta)^{{p \over k}} \\ & = & {(k!)^{{p \over k}} \over p!} \bigg({\sqrt[p]{x(x-1) \dots (x-p+1)} \over \sqrt[k]{(x)(x-1) \dots (x-k+1)}}\bigg)^p f_{k-1}(\Delta)^{{p \over k}} \\ & = & {(k!)^{{p \over k}} \over p!} \bigg(1 + {\sqrt[p]{x \dots (x-p+1)} - \sqrt[k]{(x) \dots (x-k+1)} \over \sqrt[k]{(x)(x-1) \dots (x-k+1)}}\bigg)^p f_{k-1}(\Delta)^{{p \over k}} \\ & > & {(k!)^{{p \over k}} \over p!} \bigg(1 + {\sqrt[p]{x \dots (x-p+1)} - \sqrt[p]{(x-c) \dots (x-c-p+1)} \over \sqrt[k]{(x)(x-1) \dots (x-k+1)}}\bigg)^p f_{k-1}(\Delta)^{{p \over k}} \\ & > & {(k!)^{{p \over k}} \over p!} \bigg(1 + {c \over \sqrt[k]{(x)(x-1) \dots (x-k+1)}}\bigg)^p f_{k-1}(\Delta)^{{p \over k}} \\ & = & {(k!)^{{p \over k}} \over p!} \bigg(1 + {c \over \sqrt[k]{k! f_{k-1}(\Delta)}}\bigg)^p f_{k-1}(\Delta)^{{p \over k}} \\ & = & {(k!)^{{p \over k}} \over p!} \bigg(1 + {k-p \over 2\sqrt[k]{k! f_{k-1}(\Delta)}} \bigg)^p f_{k-1}(\Delta)^{{p \over k}}.
\end{eqnarray*}
The seventh and eighth lines follow from Lemmas~\ref{narrower} and~\ref{subtractc}, respectively.  \endproof

\smallskip\noindent \textit{Proof of Corollary~\ref{noreasy}: \ }  If $f_{k-1}(\Delta) > 0$, then this follows from Theorem~\ref{withoutr} by dropping the central term, which is greater than 1.  If $f_{k-1}(\Delta) = 0$, then the left side of the inequality is positive and the right side is zero.  \endproof

Alternatively, Corollary~\ref{noreasy} follows from Theorem~\ref{withr} by taking the limit as $r \to \infty$.

\smallskip\noindent \textit{Proof of Corollary~\ref{symmetric}: \ }  Multiply the inequality of Corollary~\ref{noreasy} by $p!$ and take both sides to the ${1 \over p}$ power to get $(p! f_{p-1}(\Delta))^{{1 \over p}} > (k! f_{k-1}(\Delta))^{{1 \over k}}$.  Set $p = k-1$ and chain the inequalities to get the corollary.  \endproof

The great advantage of Theorem~\ref{withoutr} over Theorems~\ref{kktheorem} and~\ref{continuouskk} is that it is easy to compute.  There is no need to find a bunch of binomial coefficients as in Theorem~\ref{kktheorem}, nor to solve a polynomial of high degree as in Theorem~\ref{continuouskk}.  Rather, all of the computations can be done easily with a hand held calculator.  This theorem also give a far more intuitive idea of how quickly $f_{p-1}(\Delta)$ must grow as $f_{k-1}(\Delta)$ does.

The disadvantage of Theorem~\ref{withoutr} is that it isn't as sharp of a bound as Theorem~\ref{continuouskk}, let alone the sharp bounds of Theorem~\ref{kktheorem}.  Still, the bounds of Theorem~\ref{continuouskk} are usually numerically closer to those of Theorem~\ref{withoutr} than to Theorem~\ref{kktheorem}, so this is a pretty good approximation.

Corollary~\ref{noreasy} is a much weaker bound.  One virtue of the corollary is in making it intuitively obvious how quickly $f_{p-1}(\Delta)$ must grow as $f_{k-1}(\Delta)$ does.  Still, as $f_{k-1}(\Delta) \to \infty$, the ratio of the bounds of Corollary~\ref{noreasy} and the sharp bounds of the Kruskal-Katona theorem converges to 1.  Corollary~\ref{noreasy} also leads quickly to the elegant inequalities of Corollary~\ref{symmetric}.

\section{Stronger bounds}

While Theorem~\ref{withoutr} is easy to compute, it isn't that sharp of a bound.  Corollaries~\ref{noreasy} and~\ref{symmetric} are much weaker yet.  Thus, there is room for stronger bounds that are still relatively easy to compute. Theorem~\ref{withr} is just such a bound, and we prove it in this section.

To prove it, we need the notion of a colored complex.  The Frankl-F\"{u}redi-Kalai theorem characterizes the face vectors of colored complexes in much the same way that the Kruskal-Katona theorem characterizes the face vectors of simplicial complexes.  Theorem~\ref{colorapprox} is a numerical approximation to the Frankl-F\"{u}redi-Kalai theorem in a similar sense to how Theorems~\ref{continuouskk} and~\ref{withoutr} are numerical approximations to the Kruskal-Katona theorem.

Theorem~\ref{colorapprox} was proven by Frankl, F\"{u}redi, and Kalai in the same paper in which they characterized the face vectors of colored complexes.  We give an easier proof of the same theorem, written differently from how they wrote it.  We then use this approximation and the Kruskal-Katona theorem to prove Theorem~\ref{withr}.

A \textit{coloring} of a simplicial complex with color set $[r] = \{1, 2, \dots, r\}$ is an assignment of a color to each vertex of the complex such that no two vertices in the same face are the same color.  This is equivalent to no two vertices in the same edge being the same color, which is the same requirement as for a graph coloring of the 1-skeleton of the complex, taken as a graph.  If a simplicial complex can be colored with $r$ colors, then we call it \textit{$r$-colorable}.  Note that this doesn't necessarily mean that $r$ is the fewest colors possible; a complex that is $r$-colorable is automatically $(r+1)$-colorable, as we do not have to use all of the colors.

One can ask which integer vectors can arise as the face vectors of $r$-colorable complexes.  The Frankl-F\"{u}redi-Kalai theorem \cite{balanced} answers this question.  As with the Kruskal-Katona theorem, the statement of the theorem requires some background.

\begin{definition}
\textup{The Tur\'{a}n graph $T_{n,r}$ is the graph obtained by partitioning $n$ vertices into $r$ parts as evenly as possible, and making two vertices adjacent exactly if they are not in the same part.  We define ${n \choose k}_r$ to be the number of cliques on $k$ vertices in the Tur\'{a}n graph $T_{n,r}$.}
\end{definition}

\begin{lemma} \label{colorcan}
Given positive integers m, k, and r with $r \geq k$, there are unique integers $s \geq 0$ and $n_k, n_{k-1}, \dots, n_{k-s}$ such that
$$m = {n_k\choose k}_r + {n_{k-1}\choose k-1}_{r-1} + \dots + {n_{k-s}\choose k-s}_{r-s},$$
$n_{k-i}-\big\lfloor{n_{k-i}\over r-i}\big\rfloor > n_{k-i-1}$ for all $0\leq i < s,$ and $n_{k-s}\geq k-s > 0$.
\end{lemma}

\begin{theorem}[Frankl-F\"{u}redi-Kalai] \label{coloredkk}
Let $\Delta$ be an $r$-colorable simplicial complex and let $m = f_{k-1}(\Delta)$.  Let
$$m = {n_k \choose k}_r + {n_{k-1} \choose k-1}_{r-1} + \dots + {n_{k-s} \choose k-s}_{r-s}$$
as in Lemma~\ref{colorcan}.  Then
$$f_{k-2}(\Delta) \geq {n_k \choose k-1}_r + {n_{k-1} \choose k-2}_{r-1} + \dots + {n_{k-s} \choose k-s-1}_{r-s}.$$
Furthermore, if a positive integer vector with 1 as its first entry satisfies these inequalities for all $k \geq 2$, then it is the face vector of some $r$-colored simplicial complex.
\end{theorem}

The statement of this theorem is very analogous to that of the Kruskal-Katona theorem.  As before, one can program a computer to give the precise bounds.  Also like with Kruskal-Katona, these bounds can be chained, so that for $i \geq 2$,
$$f_{k-i}(\Delta) \geq {n_k \choose k-i+1}_r + {n_{k-1} \choose k-i}_{r-1} + \dots + {n_{k-s} \choose k-i-s+1}_{r-s}.$$

Unfortunately, actually doing computations with these coefficients is more awkward than with those of the Kruskal-Katona theorem.  As with Lemma~\ref{kklemma}, one can compute the constant $n_k$ by finding the unique value of $n_k$ such that ${n_k \choose k}_r \leq m < {n_k + 1 \choose k}_r$.  One then repeats the procedure with $m - {n_k \choose k}_r$ and $k-1$, and so forth.

In this setting, however, the ${n \choose k}_r$ coefficients are no longer the common binomial coefficients.  There are a variety of formulas to compute them, but even the simplest involve either summations or recursions.  One such formula is to set $p = \big\lfloor {n \over r} \big\rfloor$ and $q = n - pr$.  We can then compute
$${n \choose k}_r = \sum_{i=0}^{q} {q \choose i}{r-i \choose k-i}p^{k-i}.$$

Naturally, this extra complication makes computing the bounds of the Frankl-F\"{u}redi-Kalai theorem more awkward and less intuitive than those of the Kruskal-Katona theorem.  One can try to define ${x \choose k}_r$ to interpolate between the integer values of $x$ for a theorem analogous to Theorem~\ref{continuouskk}, but here, it isn't even immediately obvious what the formula should be.

Theorem~\ref{colorapprox} avoids the computations from Tur\'{a}n graphs entirely.  Furthermore, $r$ is given to us, so there is no messiness involved in trying to compute $r$.  Again, $r$, $p$, and $k$ are all honest integers, and ${r \choose p}$ and ${r \choose k}$ are genuine binomial coefficients.  We need a lemma before we can prove Theorem~\ref{colorapprox}.

Before the lemma, we note that ${n \choose k}_r$ is easy to compute when $n = pr$.  In this case, there are $r$ colors of vertices and $p$ vertices of each color.  To pick $k$ vertices of different colors, there are ${r \choose k}$ ways to pick the colors that will be used, and $p$ ways to pick the vertex of a given color.  Hence, ${pr \choose k}_r = {r \choose k}p^k$.

\begin{lemma} \label{colorkklimit}
Let $r \geq k > p \geq 1$ be integers and let $g(x)$ be a function whose domain is the natural numbers such that for positive integers $m$ and $n$, if ${nr \choose k}_r \leq m \leq {(n+1)r \choose k}_r$, then ${nr \choose p}_r \leq g(m) \leq {(n+1)r \choose p}_r$.  Then
$$\lim_{m \to \infty} {g(m) \over m^{p \over k}} = {r \choose p}{r \choose k}^{-{p \over k}}.$$
\end{lemma}

\proof  For each integer $m$, let $n$ be the unique integer such that ${n r \choose k}_r \leq m < {(n + 1)r \choose k}_r$.  We can compute
\begin{eqnarray*}
{{n r \choose p}_r \over {(n+1)r \choose k}_r^{{p \over k}}} \leq & {g(m) \over m^{{p \over k}}} & \leq {{(n+1) r \choose p}_r \over {n r \choose k}_r^{{p \over k}}} \\ {{r \choose p}n^p \over \big({r \choose k} (n+1)^k\big)^{{p \over k}}} \leq & {g(m) \over m^{{p \over k}}} & \leq {{r \choose p}(n + 1)^p \over \big({r \choose k} n^k\big)^{{p \over k}}} \\ {r \choose p}{r \choose k}^{-{p \over k}} \bigg({n \over n + 1}\bigg)^p \leq & {g(m) \over m^{{p \over k}}} & \leq {r \choose p}{r \choose k}^{-{p \over k}} \bigg({n + 1 \over n}\bigg)^p.
\end{eqnarray*}
As $m \to \infty$, we get $n \to \infty$ also.  As $n \to \infty$, we have both ${n \over n + 1} \to 1$ and ${n + 1 \over n} \to 1$.  Thus, both ends of the last big inequality go to ${r \choose p}{r \choose k}^{-{p \over k}}$ as $m \to \infty$.  By the squeeze theorem, so does ${g(m) \over m^{{p \over k}}}$.  \endproof

The use of this lemma is that if we set
$$m = {n_k\choose k}_r + {n_{k-1}\choose k-1}_{r-1} + \dots + {n_{k-s}\choose k-s}_{r-s}$$
as in Lemma~\ref{colorcan} and let
$$g(m) = {n_k\choose p}_r + {n_{k-1}\choose p-1}_{r-1} + \dots + {n_{k-s}\choose p-s}_{r-s},$$
then $g(x)$ satisfies the conditions of Lemma~\ref{colorkklimit}.  With this definition, Theorem~\ref{coloredkk} states that $f_{p-1}(\Delta) \geq g(f_{k-1}(\Delta))$.

We want a notion of replicating vertices to create a new complex.

\begin{definition}
\textup{Let $\Delta$ be an $r$-colored simplicial complex and let $q \geq 1$ be an integer.  Let the vertices of $\Delta$ be $v_1, v_2, \dots, v_n$.  Define a complex $\Delta^q$ on vertices $v_i^j$ for $1 \leq i \leq n$ and $1 \leq j \leq q$ such that $\{v_{i_1}^{j_1}, v_{i_2}^{j_2}, \dots, v_{i_p}^{j_p}\}$ is a face of $\Delta^q$ if and only if $\{v_{i_1}, v_{i_2}, \dots, v_{i_p}\}$ is a face of $\Delta$.}
\end{definition}

Note that this definition requires that $i_j \not = i_k$ for all $j \not = k$ in order for a vertex set to correspond to a face in $\Delta^q$.  Thus, we can make $\Delta^q$ an $r$-colored complex by giving $v_i^j$ in $\Delta^q$ the same color as $v_i$ in $\Delta$.

A face on $n$ vertices in $\Delta$ corresponds to $q^n$ faces in $\Delta^q$, as for each vertex of $\Delta$, there are $q$ ways to pick which vertex of $\Delta^q$ corresponds to it.  Thus, if the face vector of $\Delta$ is $(1, c_1, c_2, \dots, c_d)$, then the face vector of $\Delta^q$ is $(1, qc_1, q^2c_2, \dots, q^dc_d)$.

\smallskip\noindent \textit{Proof of Theorem~\ref{colorapprox}: \ }  For each positive integer $m$, let $g(m)$ be the number of faces of dimension $p-1$ of the $r$-colorable complex with the fewest such faces among all $r$-colorable complexes with exactly $m$ faces of dimension $k-1$.  Lemma~\ref{colorkklimit} asserts that
$$\lim_{m \to \infty} {g(m) \over m^{p \over k}} = {r \choose p}{r \choose k}^{-{p \over k}}.$$

We can compute
$${f_{p-1}(\Delta^q) \over (f_{k-1}(\Delta^q))^{{p \over k}}} = {f_{p-1}(\Delta) q^p \over (f_{k-1}(\Delta) q^k)^{{p \over k}}} = {f_{p-1}(\Delta) \over (f_{k-1}(\Delta))^{{p \over k}}}$$
As this holds for every positive integer $q$, we get
$$\liminf_{m \to \infty} {g(m) \over m^{p \over k}} \leq {f_{p-1}(\Delta) \over (f_{k-1}(\Delta))^{{p \over k}}}$$
Thus, we have
$${f_{p-1}(\Delta) \over (f_{k-1}(\Delta))^{{p \over k}}} \geq \liminf_{m \to \infty} {g(m) \over m^{p \over k}} = \lim_{m \to \infty} {g(m) \over m^{p \over k}} = {r \choose p}{r \choose k}^{-{p \over k}},$$
as desired.  \endproof

The Kruskal-Katona theorem doesn't merely give bounds, however. It also gives a complex that attains the bounds, by specifying that the faces of a given dimension be added in the reverse-lexicographic (``rev-lex") order. To define the rev-lex order of $i$-faces of a simplicial complex on $n$ vertices, we start by labeling the vertices $1, 2, \dots$.  Let $\mathbb{N}$ be the set of natural numbers, let $A$ and $B$ be distinct subsets of $\mathbb{N}$ with $|A| = |B| = i$, and let $A \nabla B$ be the symmetric difference of $A$ and $B$.

\begin{definition}
\textup{For $A, B \subset \mathbb{N}$ with $|A| = |B|$, we say that $A$ precedes $B$ in the rev-lex order if max$(A \nabla B) \in B$, and $B$ precedes $A$ otherwise.}
\end{definition}

For example, $\{2, 3, 5\}$ precedes $\{1, 4, 5\}$, as 3 is less than 4, and $\{3, 4, 5\}$ precedes $\{1, 2, 6\}$.

\begin{definition}
\textup{The \textit{rev-lex complex on $m$ faces of dimension $k-1$} is the pure complex whose facets are the first $m$ $k$-sets in rev-lex order.}
\end{definition}

Given positive integers $m, k,$ and $p$ with $k > p$, if $\Delta$ is the rev-lex complex on $m$ faces of dimension $k-1$ and $\Gamma$ is some other complex such that $f_{k-1}(\Gamma) = m$, the Kruskal-Katona theorem states that $f_{p-1}(\Gamma) \geq f_{p-1}(\Delta)$.

\smallskip\noindent \textit{Proof of Theorem~\ref{withr}: \ }  Let $\Delta$ be the rev-lex complex on $m$ faces of dimension $k-1$.  If $m \leq {r+1 \choose k}$, then the first $m$ faces of dimension $k-1$ in the rev-lex order contain at most $r+1$ vertices.  Thus, $\Delta$ is $(r+1)$-colorable, by the trivial coloring of giving every vertex its own color.  Furthermore, among the first ${r+1 \choose k}$ faces, the ones using both vertices $r$ and $r+1$ are the very last ones.  There are ${r-1 \choose k-2}$ possible sets of $k$ vertices among the first $r+1$ that use both vertices $r$ and $r+1$, so there are ${r+1 \choose k} - {r-1 \choose k-2} = {r \choose k} + {r-1 \choose k-1}$ that do not.  Hence, if $m \leq {r \choose k} + {r-1 \choose k-1},$ then no face of $\Delta$ contains both vertices $r$ and $r+1$.  We can make these two vertices the same color, and so $\Delta$ is $r$-colorable.  Hence
$$f_{p-1}(\Delta) \geq {r \choose p}{r \choose k}^{-{p \over k}} (f_{k-1}(\Delta))^{{p \over k}}$$
by Theorem~\ref{colorapprox}.
By the Kruskal-Katona theorem, for any complex $\Gamma$ with $f_{k-1}(\Gamma) = m$, we have
$$f_{p-1}(\Gamma) \geq f_{p-1}(\Delta) \geq {r \choose p}{r \choose k}^{-{p \over k}} (f_{k-1}(\Delta))^{{p \over k}},$$
which is what we wanted to show.  \endproof

In contrast to the Kruskal-Katona theorem, Theorem~\ref{withr} only requires us to compute one binomial coefficient, rather than as many as $k$.  In contrast to Theorem~\ref{continuouskk}, the constants in Theorem~\ref{withr} are all honest integers, so the terms that look like binomial coefficients really are.  While the bound depends on $r$, picking out the best value of $r$ is not difficult.  We compute $n_k$ as in Lemma~\ref{kklemma}, and we see if $r = n_k$ works.  It will unless $n_{k-1} = n_k - 1$ and $s \geq 2$, in which case, we use $r = n_k + 1$.

Furthermore, if $r = n_k$ works, then this lemma gives a tighter bound on $f_{p-1}(\Delta)$ than Theorem~\ref{continuouskk}.  The exception is when when $f_{k-1}(\Delta) = {r \choose k}$.  In this case, the bounds coincide and are both sharp.

In the proof of Theorem~\ref{withoutr}, we noted that the bound of Theorem~\ref{continuouskk}, namely $f_{p-1}(\Delta) \geq {x \choose p}$, could also be written as
$$f_{p-1}(\Delta) \geq {x \choose p} {x \choose k}^{-{p \over k}} f_{k-1}(\Delta)^{{p \over k}}.$$
It is not difficult to show that the quantity ${x \choose p} {x \choose k}^{-{p \over k}}$ decreases as $x$ increases.  Theorem~\ref{withr} often allows us to replace $x$ by $r = \lfloor x \rfloor$, which yields a larger bound.

\section{Flag complexes}

In this section, we prove Corollaries~\ref{flagdim} and~\ref{flagnodim}.  These corollaries rely on a theorem from a previous paper of the author \cite{myfirst} and the notion of a flag complex.

A simplicial complex is a \textit{flag complex} if every minimal non-face is a two element set.  Equivalently, if there is a set of $n$ vertices such that every two of the vertices forms an edge, then the $n$ vertices must form a face.

Flag complexes are closely related to graphs.  The clique complex of a graph is a complex such that the vertices of the complex are the same as the vertices of the graph.  A set of vertices forms a face in the clique complex exactly if it forms a clique in the graph.  It is clear that if given a set of vertices in a graph, either they form a clique or else some two of them do not form an edge.  As such, the clique complex of a graph is a flag complex.  Conversely, every flag complex is the clique complex of its 1-skeleton, taken as a graph.

\begin{theorem} \label{flagcolored}
Let $\Delta$ be a flag complex of dimension $d-1$.  Then there is a $d$-colored complex $\Gamma$ such that $f(\Gamma) = f(\Delta)$.
\end{theorem}

Because $\Gamma$ is $d$-colored, it has to follow the bounds of Theorem~\ref{coloredkk}.  This last theorem says that the flag complex $\Delta$ has to follow these same bounds, as though the complex were $d$-colorable.  If $\Delta$ has chromatic number $d$, then this is not surprising, but the chromatic number of $\Delta$ could be much larger than $d$.  This says that $\Delta$ must satisfy the bounds of Theorem~\ref{coloredkk} as though it had chromatic number $d$, even if its real chromatic number is much larger.

\smallskip\noindent \textit{Proof of Corollary~\ref{flagdim}: \ }  This follows immediately from Theorems~\ref{colorapprox} and~\ref{flagcolored}.  \endproof

\smallskip\noindent \textit{Proof of Corollary~\ref{flagnodim}: \ }  $\Delta$ cannot have a face of dimension $r$, for otherwise, we would have $f_{k-1}(\Delta) \geq {r+1 \choose k}$ from the faces of dimension $k-1$ contained in the face of dimension $r$.  Hence, the dimension of $\Delta$ is at most $r-1$.  Apply Corollary~\ref{flagdim}.  \endproof

If ${r \choose k} \leq f_{k-1}(\Delta) \leq {r \choose k} + {r-1 \choose k-1}$, then Corollary~\ref{flagnodim} coincides with Theorem~\ref{withr}, of course.  However, if ${r \choose k} + {r-1 \choose k-1} < f_{k-1}(\Delta) < {r+1 \choose k}$, then this gives a stronger bound than Theorem~\ref{withr}.  In this range, it sometimes also gives a stronger bound than the Kruskal-Katona theorem, which is sharp for simplicial complexes in general but not for flag complexes.

In the last section, we noted that Theorem~\ref{withr} often allows us to replace $x$ by $r = \lfloor x \rfloor$ in the bound
$$f_{p-1}(\Delta) \geq {x \choose p} {x \choose k}^{-{p \over k}} f_{k-1}(\Delta)^{{p \over k}}.$$
Corollary~\ref{flagnodim} says that for a flag complex, we can always make exactly this replacement.  Thus, we always get a strictly stronger bound than that of Theorem~\ref{continuouskk}, except for when the theorem is already sharp.

\section{Comparing the bounds graphically}

In this section, we give two graphs that show the numerical bounds of the various theorems proven or cited in this paper.  This gives a clear idea of how strong numerically the various bounds are.

We set $k = 10$ and $p = 7$ and compare the various bounds.  The horizontal axis is the value of $f_9(\Delta)$, while the vertical axis shows the various lower bounds on $f_6(\Delta)$.  As these are lower bounds, higher is better.

\includegraphics[scale=0.5]{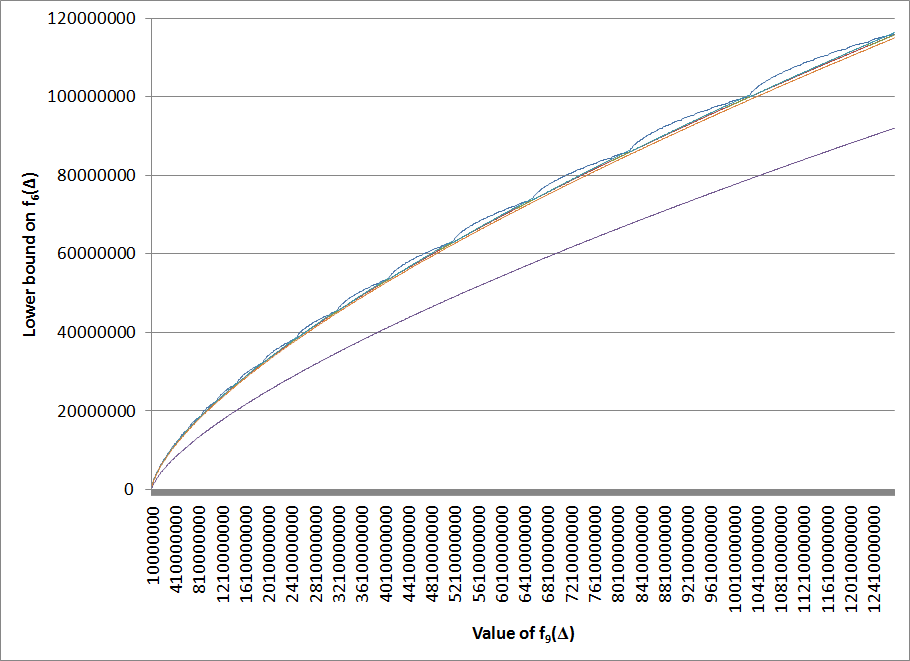}

This shows what the various bounds look like for $f_9(\Delta) < {50 \choose 10}$.  The top and somewhat jagged curve is the sharp bound of the Kruskal-Katona theorem.  The curve all the way at the bottom is that of Corollary~\ref{noreasy}.  The bounds of Theorems~\ref{withr}, \ref{withoutr}, and~\ref{continuouskk} and of Corollary~\ref{flagnodim} are also shown, but close enough together that they are hard to distinguish.  As such, we drop the bound of Corollary~\ref{noreasy} and zoom in.

\includegraphics[scale=0.5]{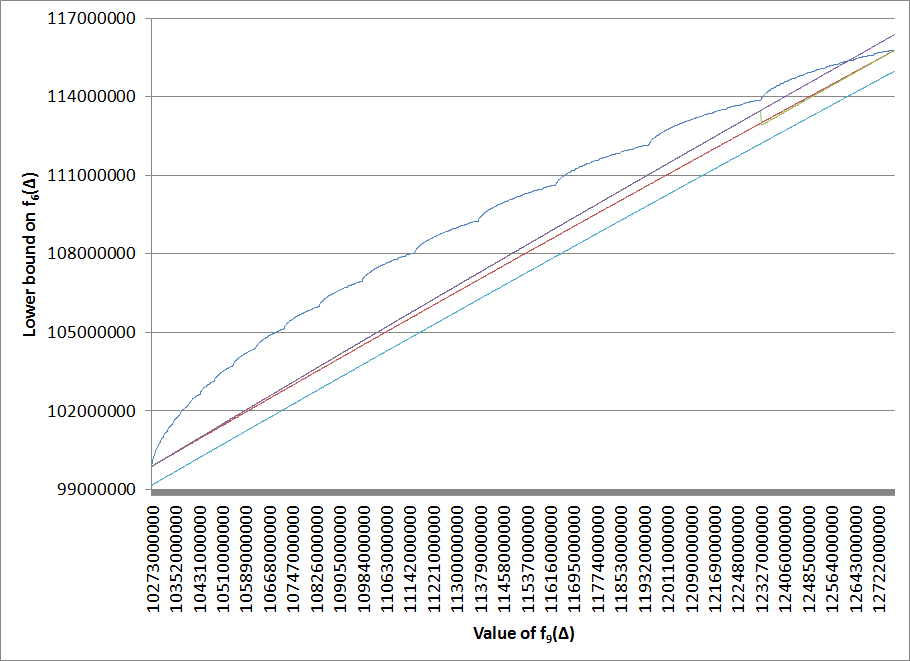}

This shows what the various bounds look like for ${50 \choose 10} < f_9(\Delta) < {51 \choose 10}$.  The jagged bound that is on top most of the way is the Kruskal-Katona theorem.  The curve at the bottom is Theorem~\ref{withoutr}.  The next lowest bound most of the way is Theorem~\ref{continuouskk}.  Theorem~\ref{withr} and Corollary~\ref{flagnodim} coincide slightly above it for ${50 \choose 10} < f_9(\Delta) \leq {50 \choose 10} + {49 \choose 9}$, at which point the bound of Theorem~\ref{withr} drops precipitously when we are forced to increase $r$ from 50 to 51.  Note that on the right side, the bound of Corollary~\ref{flagnodim} actually exceeds that of the Kruskal-Katona theorem.

\end{document}